\documentclass[11pt]{article}
\usepackage{amsmath, amsfonts, amsthm, amssymb, multicol}
\usepackage{graphicx}
\usepackage{float}
\usepackage{verbatim}

\usepackage[square,sort,comma,numbers]{natbib}
\usepackage{fancyhdr}
 
\numberwithin{equation}{section}

\newtheorem{thm}{Theorem}[section]
\newtheorem{lma}[thm]{Lemma}

\newtheorem{defn}[thm]{Definition}

\renewcommand{\geq}{\geqslant}
\renewcommand{\leq}{\leqslant}

\newcommand{\eps}{\varepsilon}

\title{ \vspace{-20mm}Almost arithmetic progressions in the primes \\ and other large sets}

\author{Jonathan M. Fraser}

\begin{document}

\date{}

\maketitle

\begin{abstract}
A celebrated and deep result of Green and Tao states that the primes contain arbitrarily long arithmetic progressions.  In this note I provide a straightforward argument demonstrating that the primes get \emph{arbitrarily close} to arbitrarily long arithmetic progressions. The argument also applies to `large sets' in the sense of Erd\H{o}s-Tur\'an.  The proof is short, completely  self-contained, and aims to give a heuristic explanation of why the primes, and other large sets,  possess arithmetic structure.
\end{abstract}

\section{Arithmetic progressions and the Green-Tao Theorem}

Arithmetic progressions are among the most natural and well-studied mathematical objects.  They are both aesthetically pleasing and ripe with structure, two properties which make it particularly interesting to find them inside other objects, which might, at first, seem  complicated and unstructured.  Such is the beauty of the Green-Tao Theorem which asserts that the primes - one of the most fundamental, complicated,  and subtle objects in mathematics  - contain \emph{arbitrarily long} arithmetic progressions.  This is a remarkable amount of additive regularity for an inherently multiplicative structure to possess.

An arithmetic progression of \emph{length} $k \geq 1$ and \emph{gap size} $\Delta>0$ is a set of the form
\[
\{x, \,  x+\Delta, \,  x+2\Delta, \,  \dots,  \, x+(k-1) \Delta \},
\]
for some $x \in \mathbb{R}$, that is, a collection of $k$ points each separated from the next by a common distance  $\Delta$. 

\begin{thm}[The Green-Tao Theorem \cite{greentao}]
The primes contain arbitrarily long arithmetic progressions, that is, for all integers $k \geq 1$ one can find an arithmetic progression of length $k$ lying somewhere in the primes.
\end{thm}

We consider the following weakened version of containing arithmetic progressions, introduced and studied in  \cite{fraseryu, frasersaitoyu}.

\begin{defn}\label{AAP}
A  set of positive integers  $X \subseteq \mathbb{Z}$  \emph{gets arbitrarily close to arbitrarily long arithmetic progressions} if, for all $k \in \mathbb{N}$ and $\varepsilon>0$, there exists an arithmetic progression $P$ of length $k$ and gap size $\Delta>0$ such that
\[
\sup_{p \in P} \inf_{x \in X} |p-x|  \leq \varepsilon \Delta.
\]
\end{defn}

This definition should be understood as saying that for arbitrarily large $k$ and arbitrarily small $\varepsilon>0$ $X$ gets within $\varepsilon$ of an arithmetic progression of length $k$.  The fact that $\varepsilon \Delta$ appears instead of $\varepsilon$ is the necessary normalization, based on the observation that all arithmetic progressions of length $k$ are essentially the same: they are all equal to $\{0,1,2, \dots, k-1\}$ upon rescaling and translation.

\begin{figure}[H]
\centering
 \includegraphics[width= 0.45\textwidth]{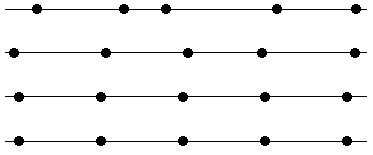}
\caption{From top row to bottom row: three different approximations to an arithmetic progression of length 5, where $\varepsilon$ is 1/3, 1/10, 1/100, respectively, followed by a genuine arithmetic progression of length $5$. At this resolution the $\varepsilon=1/100$ case is indistinguishable from the genuine arithmetic progression.}
\end{figure}

\begin{thm} \label{almostAP}
The primes get arbitrarily close to arbitrarily long arithmetic progressions.
\end{thm}

To the untrained eye this theorem may have the same aesthetic appeal as the Green-Tao Theorem: getting arbitrarily close is good enough, right? However, this theorem is very straightforward to prove, compared with the 60 page epic published in \emph{Annals of Mathematics} which is required to establish the Green-Tao Theorem \cite{greentao}.  Of course, it also follows directly from the Green-Tao Theorem.  The purpose of this article is to give a simple and self-contained proof of Theorem \ref{almostAP}.  It follows from three simple lemmas, which we will discuss in the following section and prove thereafter.

The Green-Tao Theorem is, by virtue of the fact that the sum of the reciprocals of the primes diverges (see Lemma \ref{lemma1} below),  a special case of  the Erd\H{o}s-Tur\'an conjecture on arithmetic progressions, which is a famous open problem in number theory dating back to 1936 \cite{erdos}.  It states that if $X \subseteq \mathbb{Z}^+$ is such that
\[
\sum_{x \in X} 1/x = \infty,
\]
then  $X$ should contain arbitrarily long arithmetic progressions.  We also provide a straightforward proof of the weakened version of this conjecture using the same approach.

\begin{thm}\label{erdos}
If $X \subseteq \mathbb{Z}^+$ is such that
\[
\sum_{x \in X} 1/x = \infty,
\]
then  $X$  gets arbitrarily close to arbitrarily long arithmetic progressions.
\end{thm}

This theorem was proved in \cite[Theorem 2.11]{fraseryu}, but follows directly from Lemmas \ref{lemma2} and \ref{lemma3} which we prove below.  Perhaps the interest of this result lies in the fact that, unlike in the case of the primes, the genuine version of the  Erd\H{o}s-Tur\'an conjecture is still open?

\section{Proof by three lemmas}

An obvious necessary condition for containing arbitrarily long arithmetic progressions, or even getting arbitrarily close to arbitrarily long arithmetic progressions, is being infinite.  The primes have been known to be infinite for a rather long time, the first proof often attributed to Euclid.  However, we need more in order to proceed.  Indeed, the positive  integer powers of 2 form an infinite set but it is a short exercise to see that they do \emph{not} get arbitrarily close to arbitrarily long arithmetic progressions; even arithmetic progressions of length 4!   If a sequence gets big very quickly, then the reciprocals of that sequence get small very quickly and therefore the fact that the reciprocals of the powers of 2 form a geometric series which sums to 1 is an indication of the fact that the powers of 2 grow too fast.  The following result, first proved by Euler, is a fundamental result in mathematics and is our first key ingredient.

\begin{lma}[Euler \cite{euler}] \label{lemma1}
The sum of the reciprocals of the primes diverges, that is,
\[
\sum_{p \ \textup{prime}} 1/p = \infty.
\]
\end{lma}

In the interest of being self-contained we present a well-known proof of this result due to Erd\H{o}s \cite{erdosprimes}  in Section \ref{proof1}. The next step is to turn the fact that the sum of the  reciprocals of the primes diverges  into a more quantitative statement about the distribution of the primes.  The following result is adapted from \cite[Lemma 2.10]{fraseryu} and we present a self-contained proof in Section \ref{proof2}.

\begin{lma} \label{lemma2}
If the sum of the reciprocals of a set of positive integers diverges, then the set has upper logarithmic density equal to 1, that is, if $X \subseteq \mathbb{Z}^+$ is such that
\[
\sum_{x \in X} 1/x = \infty,
\]
then 
\[
\limsup_{n \to \infty}  \sup_{m \geq 0}\frac{\log \# X \cap [m+1,m+n ]}{\log n} = 1.
\]
\end{lma}

The final step in establishing Theorem \ref{almostAP} is to show that maximal upper logarithmic density is enough to guarantee arbitrary closeness to arbitrarily long arithmetic progressions.  This result follows from \cite[Theorem 2.4]{fraseryu}, see also \cite{frasersaitoyu}, but we present a self-contained and stripped back proof in Section \ref{proof3}.

\begin{lma} \label{lemma3}
If the upper logarithmic density of a set of positive integers is equal to 1, then the set gets arbitrarily close to arbitrarily long arithmetic progressions, that is, if $X \subseteq \mathbb{Z}^+$ is such that
\[
\limsup_{n \to \infty}  \sup_{m \geq 0}\frac{\log \# X \cap [m+1,m+n ]}{\log n} = 1,
\]
then $X$ gets arbitrarily close to arbitrarily long arithmetic progressions.
\end{lma}

\section{Proofs}

\subsection{Proof of Lemma \ref{lemma1}} \label{proof1}

This proof is due to Erd\H{o}s \cite{erdosprimes} and is a classic example of proof by contradiction.  List the primes in increasing order $p_1, p_2, \dots$ and suppose that
\[
\sum_{k=1}^\infty 1/p_k < \infty
\] 
which means we can find a positive  integer $L$ such that
\[
\sum_{k=L+1}^\infty 1/p_k \leq 1/2.
\] 
Fix a large  positive integer $N$ and write $A$ for the number of integers between 1 and $N$ which are not divisible by any primes strictly larger than $p_L$.  If $n \leq N$ is such an integer then, writing $n=P^2Q$ where $Q$ is a square free integer and $P$ is an integer, one sees there are fewer than  $\sqrt{N}$ choices for $P$ and $2^L$ choices for $Q$.  Therefore
\[
A \leq \sqrt{N} 2^L.
\]
On the other hand
\[
N-A \leq \sum_{k=L+1}^\infty N/p_k  \leq N/2
\]
since there are at most $N/p_k$ positive integers less than $N$ divisible by $p_k$.  Therefore
\[
N/2 \leq A \leq \sqrt{N} 2^L
\]
which cannot be true for all $N$, yielding the desired contraction.   \qed

\subsection{Proof of Lemma \ref{lemma2}}  \label{proof2}

This proof is due to Fraser and Yu and is adapted from \cite{fraseryu}.   List the elements of $X$ in increasing order $x_1, x_2, \dots$ and suppose that the upper logarithmic density of $X$ is strictly less than 1.    It follows that there exists $s \in (0,1)$ and $C>0$ such that for all integers $m \geq 0$ and $n \geq 1$ we have
\[
 \# X \cap [m+1, m+n] \ \leq \ C n^s.
\]
For integers $N \geq 0$ write $ X_N = X\cap [2^{N}, 2^{N+1})$ and note that by the upper logarithmic density assumption
\[
\# X_N   \ \leq  \ C 2^{sN}.
\]
Therefore
\[
 \sum_{k=1}^\infty 1/x_k \ = \ \sum_{N=1}^\infty \  \sum_{k \, : \,  x_k \in X_N} 1/x_k \ \leq \  \sum_{N=1}^\infty \left( \# X_N \right) 2^{-N} \ \leq \  \sum_{N=1}^\infty C 2^{(s-1)N}   \ < \ \infty
\]
since $s<1$, which yields the desired contradiction since we assume the sum of the reciprocals of elements in $X$ diverges. \qed

\subsection{Proof of Lemma \ref{lemma3}}  \label{proof3}

 This proof is adapted from  \cite{fraseryu} and \cite{frasersaitoyu}.  Suppose $X \subseteq \mathbb{Z}^+$ does \emph{not} get arbitrarily close to arbitrarily long arithmetic progressions.  That is, there exists $k \geq 2$ and $\eps>0$ such that, given any arithmetic progression $P \subseteq \mathbb{R}$ of length $k$ and gap size $\Delta$,
\begin{equation} \label{avoid}
\sup_{p \in P} \inf_{x \in X} |p-x|  >  \eps \Delta.
\end{equation}
We may assume for convenience that $1/(2 \eps)$ is an integer, since we can always replace $\eps$ with a smaller value and force this to be true.  Fix a compact interval $J \subseteq \mathbb{R}$ of length $|J|>0$.   Cut this interval   into $k/(2 \eps)$ equal pieces of length $|J|(2 \eps)/k$ and label these from left to right by $1,2, \dots, k/(2 \eps)$.  On this set of labels (and associated intervals), form congruence classes modulo $1/(2 \eps)$ and note that the centres of the intervals with labels in the same congruence class form an arithmetic progression of length $k$ and gap size $|J|/k$.   It follows from \eqref{avoid} that at least one interval from each congruence class must not intersect $X$.  Therefore $X \cap J$ is contained in the union of  $(k-1)/(2 \eps)$ intervals of length $|J|(2 \eps)/k$.  We apply this observation inductively starting with the interval  $J_0= [m+1,m+n]$, where $m,n$ are arbitrary positive integers,  and continuing with each of the  subintervals of $J_0$ which intersect $X$.  After $N$ applications of the inductive argument we find  $X \cap J_0$ is contained in the union of  at most $\left((k-1)/(2 \eps)\right)^N$ intervals of length $(n-1)\left((2 \eps)/k\right)^N$.   Fix $N$ to be the smallest integer such that 
\[
(n-1)\left(\frac{2 \eps}{k} \right)^N < 1
\]
and note that, since $X \subseteq \mathbb{Z}$, each interval at the $N$th step contains at most one point from $X$.  It follows that 
\[
\# X \cap [m+1,  m+n]  \leq \left(\frac{k-1}{2 \eps} \right)^N \leq  \left(\frac{k-1}{2 \eps} \right)^{\frac{\log(n-1)}{\log(k/(2 \eps))}+1} \leq  \left(\frac{k-1}{2 \eps} \right) n^{\frac{\log\left((k-1)/(2 \eps)\right)}{\log(k/(2 \eps))}} 
\]
which proves that the upper logarithmic density of $X$ is bounded above by 
\[
\frac{\log\left(\frac{k-1}{2 \eps}\right)}{\log \left( \frac{k}{2 \eps} \right)} < 1
\]
contradicting our assumption that the upper logarithmic density of $X$ is equal to 1.  \qed


\vfill

\begin{centering}

\textbf{Acknowledgments}

The author was financially  supported by a \emph{Leverhulme Trust Research Fellowship} (RF-2016-500) and  an \emph{EPSRC Standard Grant} (EP/R015104/1).  He thanks Han Yu for many inspiring conversations related to the topics presented here.
\end{centering}

\vspace{5mm}

\noindent \emph{Jonathan M. Fraser\\
School of Mathematics and Statistics\\
The University of St Andrews\\
St Andrews, KY16 9SS, Scotland} \\
\noindent  Email: jmf32@st-andrews.ac.uk

\end{document}